\documentstyle[12pt]{article}

\def\be{\begin{equation}}
\def\ee{\end{equation}}
\def\bea{\begin{eqnarray}}
\def\eea{\end{eqnarray}}

\def\R{{\rm I\kern-.2em R}}
\def\1{\'{\i}}

\def\aone{a_+}

\def\athree{b_+}

\def\afive{a}
\def\asix{b}

\def\fama{${\mbox{I}_+}$}

\def\famc{${\mbox{II}}$}

\def\K{J_3}
\def\P{J_+}
\def\T{J_-}
\def\M{I}

\def\mmm{M}

\begin{document}

\thispagestyle{empty}

 \
\hfill\

\
\vspace{1cm}

\begin{center}
{\LARGE{\bf{Quantum harmonic  oscillator algebras}}}

{\LARGE{\bf{as non-relativistic
limits of}}}

{\LARGE{\bf{multiparametric $gl(2)$ quantizations}}}
\end{center}

\bigskip\bigskip

\begin{center} Angel Ballesteros$^\dagger$, Francisco J.
Herranz$^\dagger$  and Preeti Parashar$^{\dagger\ddagger}$
\end{center}

\begin{center} {\it {$^\dagger$Departamento de F\1sica, Universidad
de Burgos} \\   Pza. Misael Ba\~nuelos,
E-09001 Burgos, Spain}
\end{center}

\begin{center} {\it {$^\ddagger$Departamento de F\1sica Te\'orica,
Universidad Aut\'onoma de Madrid} \\   Cantoblanco,
E-28049 Madrid, Spain}
\end{center}

\bigskip

\begin{abstract}
Multiparametric quantum $gl(2)$ algebras are  presented
according to a classification based on their corresponding Lie
bialgebra structures. From  them, the non-relativistic limit
leading to  quantum harmonic oscillator algebras is
implemented in  the form of generalized Lie bialgebra
contractions. 
\end{abstract}

\bigskip\bigskip

\newpage


\section{Introduction}

The $gl(2)$ Lie algebra can be viewed
as the natural relativistic analogue of the
one-dimensional harmonic oscillator algebra $h_4$ \cite{Ala}.
Reciprocally, $h_4$ can be obtained from $gl(2)$  through a
generalized In\"on\"u-Wigner contraction that translates into
mathematical terms the non-relativistic limit $c\to \infty$.  
Explicitly, if we consider the commutation  relations and
second-order Casimir of the $gl(2)$ Lie algebra
\bea
&&[\K,\P]=2\P ,\quad [\K,\T]=-2\T ,\quad [\P,\T]=\K ,\quad
[\M,\cdot\,]=0 ,\cr
&&{\cal C}=\K^2 + 2 \P \T + 2 \T \P,
\label{ab}
\eea
and we apply the map defined by
\be
A_+ = \varepsilon \P ,\qquad A_- = \varepsilon \T ,\qquad
N =  (\K + \M)/2 ,\qquad  \mmm = \varepsilon^2 \M ,
\label{na}
\ee
then the limit
$\varepsilon\to 0$ ($\varepsilon=1/c$) leads to  the harmonic
oscillator algebra
$h_4$
\be
[N,A_+]=A_+ ,\qquad [N,A_-]=-A_- ,\qquad
[A_-,A_+]=\mmm ,
\qquad [\mmm,\cdot\,]=0 .
\label{nb}
\ee
The Casimir of $h_4$ is also obtained by computing 
$\lim_{\varepsilon\to 0}
\frac 12 {\varepsilon^2}  ( - {{\cal C}} 
+\M^2 )$:
\be
{\cal C}=2N\mmm -A_+A_- - A_- A_+ .
\ee

Recently, a systematic and
constructive approach to multiparametric  quantum $gl(2)$
algebras based on the classification of their associated Lie
bialgebra structures has been presented \cite{bhp}. In that paper, the
question concerning the generalization of the Lie bialgebra
contraction procedure to multiparametric structures has been
also solved. Now, we make use of those results in order to obtain
several quantum $h_4$ algebras and their associated deformed
Casimir operators. We emphasize that all these quantum
$h_4$ algebras are endowed with a Hopf algebra  structure,
which can be related to integrability
properties of associated models. In particular, note that the
quantum group symmetry of the 
spin $1/2$ Heisenberg XXZ and XXX chains with twisted periodic boundary
conditions \cite{AGR,PS} is given by quantum
$gl(2)$ algebras \cite{bhp,MRplb} whose non-relativistic limit will be 
analysed.


\section{Quantum $gl(2)$ algebras}

In this section we  present some relevant quantum 
 ${gl}(2)$ Hopf algebras \cite{bhp}.
Deformed Casimir operators, essential for the construction of
integrable systems
\cite{orl}, and quantum $R$-matrices are also explicitly given.

\subsection{Family \fama\ quantizations}

\subsubsection{Standard subfamily  $U_{\aone,\afive}({gl}(2))$
with $\aone\ne 0$, $\afive\ne 0$}

The  
quantum algebra $U_{\aone,\afive}({gl}(2))$ and its Casimir are
given by 
\bea
&&\!\!\!\!
\!\!\!\!
\Delta(\K')=1\otimes \K' + \K'\otimes 1 ,\qquad
\Delta(\P) =e^{ \afive \K' /2}\otimes \P + \P\otimes
e^{- \afive \K' /2} ,\cr
&&\!\!\!\!
\!\!\!\!
\Delta(\M)=1\otimes \M + \M\otimes 1 ,
\qquad
\Delta(\T) =e^{ \afive \K' /2}\otimes \T + \T\otimes
e^{- \afive \K' /2} ,\cr
&&\!\!\!\!
\!\!\!\!
[\K',\P]=2\P ,\quad [\K',\T]=-2\T -\frac{\aone}{\afive}\,
\frac{\sinh(\afive \K'/2)}{\afive/2}  -
\frac{\aone^2}{\afive^2}\P ,\quad [\M,\cdot\,]=0 ,\cr
&&\!\!\!\!
\!\!\!\!
  [\P,\T]=\frac{\sinh \afive \K' }{\afive}
+\frac{\aone}{ \afive} \left(\frac{e^{\afive}-1}{2\afive}\right)
\left(e^{ -\afive \K' /2}\P + \P e^{ \afive \K' /2}\right) ,
\label{ba}\\
&&\!\!\!\!
\!\!\!\!
{\cal
C}_{\aone,\afive}=
\frac{2}{\afive\tanh\afive}\left(\cosh(\afive\K')
-1\right)+ 2 (\P \T +
\T \P)+\frac{\aone^2}{\afive^2}\P^2\cr
&&\!\!\!\!
\!\!\!\!
\qquad\quad +
\frac{\aone}{\afive}
\left(\frac{\sinh(\afive \K'/2)}{\afive/2}\P
+\P\frac{\sinh(\afive \K'/2)}{\afive/2}\right) ,
\nonumber 
\eea
where $\K'=\K- \frac{\aone}{\afive}\P$.
This  quantum algebra is just
a superposition of the standard and non-standard  deformations
of $sl(2,\R)$ since  the underlying standard classical
$r$-matrix is  $r=\frac 12(\aone \K'\wedge \P - 2 \afive
\P\wedge \T)$. This fact can be clearly appreciated by
considering  the  $4\times 4$ quantum
$R$-matrix associated to
$U_{\aone,\afive}({gl}(2))$ \cite{bhp}: 
\be
{\cal R}=\left(\begin{array}{cccc}
1&h&-qh&h^2\\
0&q&1-q^2&qh\\
0&0&q&-h\\
0&0&0&1
\end{array}\right) ,\qquad 
q=e^{\afive}, \qquad h=\frac{\aone}{2}\left(\frac
{e^{\afive}-1}{\afive}\right) .
\ee
The limit $\aone\to 0$ yields the standard 
$R$-matrix of
$sl(2,\R)$, while taking
$\afive\rightarrow 0$ gives rise to the
non-standard  one. This quantum  algebra
underlies the construction of non-standard $R$-matrices out of
standard ones introduced in  \cite{AKS,ACC}.

\subsubsection{Non-standard  subfamily
$U_{\aone,\athree}({gl}(2))$ with $\aone\ne 0$}

The Hopf algebra  $U_{\aone,\athree}({gl}(2))$, whose Lie
bialgebra is generated by the triangular classical 
$r$-matrix $r=\frac
12(\aone \K\wedge \P +\athree \P\wedge \M)$, 
is given by
\bea
&&\Delta(\P)=1\otimes \P + \P\otimes 1 ,\qquad
\Delta(\M)=1\otimes \M + \M\otimes 1 ,\cr
&&\Delta(\K)=1\otimes \K + \K\otimes e^{\aone \P} -
\athree \M\otimes \left(\frac {e^{\aone \P}-1}{\aone}\right),\cr
&&\Delta(\T)=1\otimes \T + \T\otimes e^{\aone \P}
- \frac {\athree}2 \left( \K -
\frac {\athree}{\aone} \M \right)\otimes \M e^{\aone \P} ,\cr
&&[\K, \P] =2 \frac {e^{\aone \P} - 1} {\aone} , \qquad
[\K, \T] = - 2 \T + \frac{\aone}{2}\left(
\K - \frac{\athree}{\aone}\M\right)^2 , \cr
&&[\P, \T] = \K +  {\athree} \M 
 \frac { e^{\aone \P} - 1}{\aone}
, \qquad [\M,\,\cdot \,] = 0,
\label{bb}\\
&&{\cal C}_{\aone,\athree}=
\left(\K-\frac{\athree}{\aone} \M\right)
e^{-\aone\P}\left(\K-\frac{\athree}{\aone} \M\right)
+ 2 \frac{\athree}{\aone} \K\M\cr
&&\qquad
+2\frac{1-e^{-\aone\P} }{\aone}\T+
2\T \frac{1-e^{-\aone\P} }{\aone}
+2(e^{-\aone\P}-1).
\nonumber
\eea
This quantum algebra has  been also
obtained in \cite{Dobrev,boson,Preeti} and 
  its universal quantum $R$-matrix can be found in 
\cite{boson,Preeti}.

\subsection{Family \famc\ quantizations}

\subsubsection{Standard subfamily $U_{\afive,\asix}({gl}(2))$
 with $\afive\ne 0$}

The corresponding coproduct, commutation rules  and Casimir are
given by
\bea
&&\Delta(\M)=1\otimes \M + \M\otimes 1 ,\qquad
 \Delta(\K)=1\otimes \K + \K\otimes 1 ,\cr
&&\Delta(\P)=e^{(\afive \K - \asix \M)/2}\otimes \P + \P\otimes
e^{-(\afive \K - \asix \M)/2}  ,\cr
&&
 \Delta(\T)=e^{(\afive \K + \asix \M)/2}\otimes \T + \T\otimes
e^{-(\afive \K + \asix \M)/2} ,
\label{bc}\\
&&[\K,\P]=2\P ,\quad [\K,\T]=-2\T ,\quad [\P,\T]=\frac{\sinh
\afive\K}{\afive},
\quad [\M,\cdot\,]=0,\cr
&&
{\cal C}_{\afive}=\cosh \afive \left(\frac{\sinh
(\afive\K/2)}{\afive/2} \right)^2 
 +2\,\frac{\sinh \afive}{\afive}\,(\P
\T + \T \P).
\nonumber
\eea
This quantum algebra, together with its  universal quantum
$R$-matrix, has been obtained in \cite{CJ};  it is just the
quantum algebra underlying the XXZ Heisenberg Hamiltonian with
twisted boundary conditions \cite{MRplb}. This
deformation can be thought of as a Reshetikhin twist of the usual
standard deformation since 
in the associated   $r$-matrix, $r=- \frac 12  \asix
\K\wedge \M - \afive \P\wedge\T$,  the second
term generates the standard deformation and the  exponential of
the first one gives us the Reshetikhin twist.

\subsubsection{Non-standard subfamily 
$U_{\athree,\asix}({gl}(2))$}

The   coproduct reads
\bea
&&\Delta(\M)=1\otimes \M + \M \otimes 1 ,\qquad
\Delta(\P)=1\otimes \P + \P \otimes e^{\asix \M} ,\cr
&&\Delta(\K)=1\otimes \K + \K \otimes 1 + 
\athree \P\otimes \left(\frac
{e^{\asix \M}-1}{\asix}\right) ,\cr
&&\Delta(\T)=1\otimes \T + \T \otimes e^{-\asix \M} 
+ \athree \K\otimes
\left(\frac {e^{-\asix \M}-1}{2\asix}\right) \cr
&&\qquad \qquad+ \athree^2 \P\otimes
\left(\frac {1 - \cosh {\asix \M}}{2\asix^2}\right) ,
\label{bd}
\eea
and the associated commutation rules and Casimir are
non-deformed ones (\ref{ab}). A twisted XXX Heisenberg Hamiltonian invariant
under 
$U_{\athree,\asix}({gl}(2))$ has been constructed in \cite{bhp}.
The $r$-matrix is  $r=-\frac 12
(\asix \K - \athree
\P)\wedge \M$ and the  universal $R$-matrix turns out to be
${\cal R}= \exp\{r\}$,  which in the fundamental representation 
reads
\be
{\cal R}=\left(\begin{array}{cccc}
1&-e^{-\asix}\,p&p&-e^{-\asix}\,p^2\\
0&e^{-\asix}&0&e^{-\asix}\,p\\
0&0&e^{\asix}&-p\\
0&0&0&1
\end{array}\right)  ,\qquad
p=\frac{\athree}{2}\left(\frac{e^{\asix}-1}{ \asix}\right).
\ee


\section{Contractions to quantum oscillator algebras}

In the sequel, we work out the contractions from the above quantum $gl(2)$ 
algebras to  quantum $h_4$ algebras (a systematic
approach to the latter structures can be found in \cite{osc}).  In order to
contract a given quantum algebra we have to consider  the
In\"on\"u-Wigner contraction (e.g.\ (\ref{na})) together with
a   mapping $a=\varepsilon^{n} a'$ on {\em each} initial deformation parameter
$a$
 where $n$ is any real number and $a'$ is the
contracted deformation parameter  \cite{LBC}. The convergency  of both
the classical
$r$-matrix and the cocommutator $\delta$ under the limit
$\varepsilon \to 0$ have to be analysed
separately, since  starting 
from a coboundary bialgebra, the contraction  can lead  to either
another coboundary bialgebra (both $r$ and $\delta$ converge)
or to  a non-coboundary one  ($r$ diverges but
$\delta$ converges).  Hence we have to find out the
minimal value of the number $n$  such that $r$ converges, 
the minimal value of $n$ such that $\delta$ converges,
and finally to compare both of them.

\subsection  {Standard family \famc: 
$U_{\afive,\asix}({gl}(2))\to U_{\xi,\vartheta}(h_4)$ with
$\xi\ne 0$}
 
Let us illustrate our procedure starting with 
the quantum algebra $U_{\afive,\asix}({gl}(2))$.
We consider the
maps
\be
\afive= -\varepsilon^{n_{\afive}}\xi ,
\qquad \asix=-\varepsilon^{n_{\asix}}\vartheta  ,
\label{sa}
\ee
where $\vartheta$, $\xi$ are the contracted deformation
parameters, and $n_{\afive}$, $n_{\asix}$ are real numbers to
be determined by imposing  the convergency of $r$. We introduce
the maps (\ref{na}) and (\ref{sa}) in the  classical $r$-matrix
associated to
$U_{\afive,\asix}({gl}(2))$:
\be 
\begin{array}{l}
 r=- \frac 12  \asix \K\wedge \M  - \afive \P\wedge\T\cr
 \quad =
\frac 12 \varepsilon^{n_{\asix}}\vartheta 
(2 N - \mmm  \varepsilon^{-2})\wedge \mmm  \varepsilon^{-2} 
+\varepsilon^{n_{\afive}}\xi
  A_+ \varepsilon^{-1}\wedge A_- \varepsilon^{-1}\cr
 \quad =  \varepsilon^{n_{\asix}-2}\vartheta 
  N  \wedge \mmm   
+\varepsilon^{n_{\afive}-2}\xi 
  A_+ \wedge A_-  .
\end{array}
\ee 
Hence the minimal values of the indices $n_{\afive}$,
$n_{\asix}$ which ensure the convergency of $r$ under the limit
$\varepsilon \to 0$ are $n_{\afive}=2$,
$n_{\asix}=2$. Now we have to analyse the convergency of the
cocommutator $\delta$ associated to
$U_{\afive,\asix}({gl}(2))$. Thus we consider  the maps 
(\ref{sa}) and look for the minimal values of 
$n_{\afive}$, $n_{\asix}$ which allow $\delta$ to converge
under the limit $\varepsilon \to 0$. It can be checked that
they are again $n_{\afive}=2$,
$n_{\asix}=2$, so that the resulting $h_4$   bialgebra is  
coboundary  
(both $n_{\afive}$, $n_{\asix}$ coincide for $r$ and $\delta$).
Therefore  
the  transformations of the  deformation parameters 
so obtained are
$\afive= -\varepsilon^2\xi$ and
$\asix=-\varepsilon^2\vartheta$. 
Finally, we introduce these maps together with (\ref{na}) in  
 $U_{\afive,\asix}({gl}(2))$ and we obtain the
following quantum oscillator algebra 
$U_{\xi,\vartheta}(h_4)$:
\bea
&&\Delta(N)=1\otimes N + N\otimes 1 ,\qquad
\Delta(\mmm)=1\otimes \mmm + \mmm\otimes 1 ,\cr
&&\Delta(A_+)=e^{(\vartheta+\xi)\mmm/2}
\otimes A_+ + A_+\otimes
e^{-(\vartheta+\xi)\mmm/2} ,\cr
&&\Delta(A_-)=e^{-(\vartheta-\xi)\mmm/2}\otimes A_- + A_-\otimes 
e^{(\vartheta-\xi)\mmm/2},
\eea
\be
[N,A_+]=A_+,\quad
[N,A_-]=-A_-,\quad 
[A_-,A_+]=\frac{\sinh \xi \mmm}{\xi},\quad
[\mmm,\cdot\,]=0 .
\ee
The deformed oscillator Casimir comes from $\lim_{\varepsilon\to 0}
\frac 12{\varepsilon^2}  \bigl( -{\cal C}_{\afive} +
  \bigl(\frac{\sinh (\afive \M/2)}{\afive/2}\bigr)^2 \bigr)$:
\be
{\cal C}_{\xi}=2 N \frac{\sinh \xi \mmm}{\xi}- A_+A_- - A_-A_+ .
\ee 
If $\vartheta=0$, the quantum oscillator introduced
in
\cite{Enrico,GS} is recovered.

Hereafter we give the transformations of the deformation 
parameters for the remaining quantum $gl(2)$ algebras together
with the resulting quantum $h_4$ algebras; we stress that in all
cases the contractions are found to have a coboundary character.

\subsection {Non-standard family \famc: 
$U_{\athree,\asix}({gl}(2))\to U_{\beta_+,\vartheta}(h_4)$}

The transformations of the deformation parameters are
$\athree= 2\varepsilon^3\beta_+$ and
$\asix=-\varepsilon^2\vartheta$. The coproduct of the quantum
oscillator algebra
$U_{\beta_+,\vartheta}(h_4)$ reads
\bea
&&
\Delta(\mmm)=1\otimes \mmm + \mmm\otimes 1,\qquad
\Delta(A_+)=1\otimes A_+ + A_+\otimes
e^{- \vartheta \mmm},\cr
&&\Delta(A_-)=1\otimes A_- + A_-\otimes  e^{\vartheta \mmm}
+\beta_+ \mmm \otimes \left(\frac{e^{\vartheta \mmm} - 1}{\vartheta} 
\right),\cr
&&\Delta(N)=1\otimes N + N\otimes 1+ \beta_+ A_+\otimes
\left(\frac{1-  e^{-\vartheta \mmm}}{\vartheta} 
\right) .
\eea
Commutation rules and Casimir of 
$U_{\beta_+,\vartheta}(h_4)$ are the  non-deformed
ones  (\ref{nb}).

\subsection  {Standard family \fama: 
$U_{\aone,\afive}({gl}(2))\to U_{\beta_+,\xi}(h_4)\to 
U_{\xi}(h_4)$ with $\xi\ne 0$}

In this case,   the maps  
$\aone = 2 \varepsilon^3 \beta_+$ and $\afive=-\varepsilon^2 
\xi$  lead to  
$U_{\beta_+,\xi}(h_4)$:
\bea
&& 
\!\!\!\!
\!\!\!\!
\Delta(\mmm)=1\otimes \mmm + \mmm\otimes 1 ,
\qquad \Delta(A_\pm)=e^{ \xi\mmm/2}
\otimes A_\pm + A_\pm \otimes
e^{- \xi\mmm/2} ,\cr
&&\!\!\!\!
\!\!\!\!
\Delta(N)=1\otimes N + N \otimes 1
+\beta_+\left(\frac{1-e^{\xi\mmm/2}}{\xi}\right)\otimes A_+
+\beta_+ A_+\otimes
\left(\frac{1-e^{-\xi\mmm/2}}{\xi}\right),\cr 
&&\!\!\!\!
\!\!\!\!
[N,A_+]=A_+,\qquad  [A_-,A_+]=\frac{\sinh \xi
\mmm}{\xi},\qquad [\mmm,\cdot\,]=0 ,
\label{masi}\\ 
&&\!\!\!\!
\!\!\!\!
[N,A_-]=-A_- + \frac{\beta_+}{\xi}\left(
\frac{\sinh \xi \mmm}{\xi} - \frac{\sinh (\xi \mmm/2)}{\xi/2}
\right),\cr 
&&\!\!\!\!
\!\!\!\!
{\cal C}_{\beta_+,\xi}=2 N \frac{\sinh \xi \mmm}{\xi}
- A_+A_- - A_-A_+ +
2 A_+ \frac{\beta_+}{\xi}\left(
\frac{\sinh \xi \mmm}{\xi} - \frac{\sinh (\xi \mmm/2)}{\xi/2}
\right) ,
\nonumber
\eea
where the Casimir is provided by $\lim_{\varepsilon\to 0}
\frac 12{\varepsilon^2}  \bigl( -{\cal C}_{\aone, \afive} +
  \bigl(\frac{\sinh (\afive \M/2)}{\afive/2}\bigr)^2 \bigr)$.
 However
the parameter $\beta_+$ is irrelevant and it can be removed
from (\ref{masi})  by applying the change of basis defined by
\be
N'=N+\frac{\beta_+}{\xi}A_+ ,\quad A'_+=A_+ ,
\quad
A'_-=A_-+\frac{\beta_+}{\xi}\frac{\sinh(\xi\mmm/2)}{\xi/2}
,\quad
\mmm'=\mmm .
\ee
Thus we recover  $U_{\xi}(h_4)$, already obtained in
sec.\ 3.1  as
$U_{\vartheta,\xi}(h_4)
\to U_{\vartheta=0,\xi}(h_4)$.

\subsection {Non-standard family \fama: 
$U_{\aone,\athree}({gl}(2))\to U_{\alpha_+}(h_4)$ with \break
$\alpha_+\ne 0$}

The transformations of the deformation  parameters turn out to
be  $\aone=\varepsilon \alpha_+$ and 
$\athree=-\varepsilon\alpha_+$. 
Hence, we obtain the  ``Jordanian $q$-oscillator" \cite{osc} 
$U_{\alpha_+}(h_4)$:
\bea
&&\Delta(A_+)=1\otimes A_+ + A_+\otimes 1,\qquad
\Delta(\mmm)=1\otimes \mmm + \mmm\otimes 1 ,\cr
&&\Delta(A_-)=1\otimes A_- + A_-\otimes e^{\alpha_+ A_+}
+\alpha_+ N\otimes \mmm e^{\alpha_+ A_+},\cr
&&\Delta(N)=1\otimes N+ N\otimes e^{\alpha_+ A_+},
\qquad [\mmm,\cdot\,]=0,\\
&&[N,A_+]=\frac{e^{\alpha_+ A_+}-1}{\alpha_+},\quad
[N,A_-]=-A_-,\quad 
[A_-,A_+]=\mmm e^{\alpha_+ A_+}.
\nonumber
\eea
The quantum Casimir is computed as 
$\lim_{\varepsilon\to 0}
\varepsilon^2\left( - \frac{1}{2}{{\cal C}_{\aone,\athree}} 
+\M^2\right)$, 
and reads,
\be
{\cal C}_{\alpha_+}=2N\mmm +\frac{e^{-\alpha_+ A_+ }- 1}{\alpha_+} A_-
+A_- \frac{e^{-\alpha_+ A_+} - 1}{\alpha_+} .
\ee


\bigskip
\bigskip

\noindent
{\Large{{\bf Acknowledgments}}}

\bigskip

 A.B. and
F.J.H. have been partially supported by DGICYT (Project  PB94--1115)
from the Ministerio de Educaci\'on y Cultura  de Espa\~na and by Junta
de Castilla y Le\'on (Projects CO1/396 and CO2/297). P.P. has been
supported by a fellowship from AECI, Spain.


\bigskip
\bigskip





\begin{thebibliography}{99}

\bibitem{Ala}  Aldaya V.,  Bisquert J., and  Navarro-Salas J.: 
  Phys. Lett. A {\em 156} (1991) 381.


\bibitem{bhp}  Ballesteros A., Herranz F.\,J., and Parashar P.:
preprint math.QA/9806149.


\bibitem{AGR}  Alcaraz F., Grimm  U., and Rittenberg  V.:  
 Nucl. Phys. B {\em 316} (1989) 735.

\bibitem{PS} 
 Pasquier V. and  Saleur H.:  Nucl. Phys. B {\em 330} (1990) 523.



\bibitem{MRplb}   Monteiro M.\,R.,  Roditi I.,  Rodrigues
L.\,M.\,C.\,S., and  Sciuto S.:
  Phys. Lett. B {\em 354} (1995) 389.



\bibitem{orl}  Ballesteros A. and  Ragnisco O.: 
  J. Phys. A: Math. Gen.  {\em 31} (1998) 3791.








\bibitem{AKS}
 Aghamohammadi A., Khorrami  M., and  Shariati  A.:
J. Phys. A: Math. Gen.  {\em 28} (1995) L225.

\bibitem{ACC}
Abdesselam B.,  Chakrabarti A., and  Chakrabarti R.:
  Mod. Phys. Lett. A {\em 13} (1998) 779.

\bibitem{Dobrev}  Aneva B.L.,  Dobrev V.K., and Mihov  S.G.: 
 J. Phys. A: Math. Gen.  {\em 30} (1997) 6769.

\bibitem{boson}   Ballesteros A.,  Herranz F.J., and  Negro J.:
 J. Phys. A: Math. Gen.  {\em 30} (1997) 6797.



\bibitem{Preeti}
     Parashar P.:   Lett. Math. Phys. to appear,
preprint q-alg/9705027.

\bibitem{CJ} Chakrabarti R.  and  Jagannathan R.:  J. Phys. A:
Math. Gen.  {\em 27} (1994) 2023.

\bibitem{osc}
 Ballesteros A. and    Herranz F.J.:   J. Phys. A: Math.
Gen.  {\em 29} (1996) 4307.


\bibitem{LBC}
 Ballesteros  A.,  Gromov N.A.,  Herranz F.J.,   del Olmo
M.A.,  and  Santander  M.:   J. Math. Phys.   {\em 36}
(1995)  5916.

\bibitem{Enrico}
 Celeghini E., Giachetti R., Sorace E. and Tarlini M.: J. Math. Phys. {\em
36} (1991)  1155.

\bibitem{GS}
 G\'omez C. and Sierra G.:   J. Math. Phys.   {\em 34}
(1993)  2119.





\end{thebibliography}
\end{document}